\documentclass{amsart}[12]
\usepackage{mathrsfs, amsmath,amsfonts,amssymb}
\usepackage{amsmath}
\usepackage{mathtools}
\usepackage{graphicx}
\usepackage[left=3cm, right=2.5cm, top=3.5cm, bottom=2cm]{geometry}
\newtheorem{thm}{Theorem}[subsection]
\newtheorem{prop}[thm]{Proposition}
\newtheorem{lem}[thm]{Lemma}

\newtheorem{rem}[thm]{Remark}

\newcommand{\E}{\mathbb{E}}
\newcommand{\C}{\mathbb{C}}
\renewcommand{\P}{\mathbb{P}}
\newcommand{\R}{\mathbb{R}}
\newcommand{\Z}{\mathbb{Z}}
\newcommand{\dilog}{\text{dilog}}
\newcommand{\K}{\mathcal{K}}
\usepackage{fullpage}
\title{A law of large number for local patterns in plane partitions}

\author[P. Lazag]{Pierre Lazag}
\address{I2M, CNRS, Aix-Marseille Universit\'e \\Technop\^ole de Ch\^ateau-Gombert \\ 13453 Marseille cedex 13 \\ France}
\email{pierre.lazag@univ-amu.fr}

\begin{document}
\begin{abstract}
In this note, we prove a law of large numbers for local patterns in plane partitions with geometric weight, a model introcued by Okounkov and Reshetikhin. Its proof is based on the determinantal structure of the corresponding point process, wich allows to control the decay of the correlations in a convenient manner.
\end{abstract}

\maketitle

\section{Introduction}
\subsection{Introduction of the model}
We denote by $\mathbb{N}$ the set of positive integers, and by $\mathbb{N}_0$ the set of non-negative integers. A plane partition is a double sequence of non-increasing integers whith a finite number of non-zero elements. More precisely, $\pi=(\pi_{i,j})_{i,j \geq 1} \in \mathbb{N}_0^{\mathbb{N}\times \mathbb{N}}$ is a plane partition if and only if :
\begin{align*}
\pi_{i+1,j} & \leq \pi_{i,j} \text{ and } \pi_{i,j+1} \leq \pi_{i,j} \text{ for all } i,j \in \mathbb{N}, \\
|\pi| & := \sum_{i,j \geq 1} \pi_{i,j} < +\infty.
\end{align*}
To a plane partition we associate a subset of $E:=\Z\times \frac{1}{2}\Z$ via the map :
\begin{align*}
\pi \mapsto \mathfrak{S}(\pi):= \{ (i-j,\pi_{i,j}-(i+j-1)/2), \hspace{0.1cm} i,j \geq 1 \}.
\end{align*}
Our probability space is the space of configuration on $E$ : $\Omega := \{0,1\}^E$
equipped with the usual Borel structure generated by the cylinders. The first coordinate of a point $(t,h) \in E$ might be interpreted as the time coordinate, and the second as the space coordinate. For a plane partition $\pi$ and $(t,h) \in E$, we define $c_{(t,h)}(\pi) \in \Omega$ by :
\begin{align*}
c_{(t,h)}(\pi)=1 \text{ if and only if } (t,h) \in \mathfrak{S}(\pi),
\end{align*}
and for a subset $m= \{ m_1,...,m_l \} \subset E$, the configuration $c_m(\pi)$ is the product :
\begin{align*}
c_m=c_{m_1}...c_{m_l}.
\end{align*}
For $q \in (0,1)$, we consider the geometric probability measure $\mathbb{P}_q$ on the set of all plane partition given by:
\begin{align*}
\mathbb{P}_q(\pi) = M q^{|\pi|}
\end{align*}
where $M$ is the normalization constant given by MacMahon formula (\cite{stanley}, corollary 7.20.3) :
\begin{align*}
M=\prod_{n=1}^{+\infty}(1-q^n)^n.
\end{align*}
By the inclusion-exclusion principle, $\P_q$ is fully characterized by the quantities :
\begin{align*}
\E_q \left[ c_m \right]= \P_q \left( m \subset \mathfrak{S}(\pi) \right), \quad m \subset E \text{, $m$ is finite}. 
\end{align*}
The main result of this note, Theorem \ref{thm1} below, is a law of large numbers, as $q=e^{-r}$ tends to $1$, for local patterns $m$ in plane partitions distributed according to $\P_q$ ; it states that the normalized sum :
\begin{align*}
\Sigma(f,m,r) = r^2 \sum_{(t,h)} f(rt,rh) c_{(t,h)+m},
\end{align*}
converges with respect to $\P_q$ to a constant, as $q=e^{-r}$ tends to $1$. Here the sum is taken over a subset of $E$ which can be regarded as the limit shape of a large plane partition scaled by a factor $1/r^2$, and $f$ is a continuous compactly supported function defined on the plane. The constant is explicitly given in terms of the discrete extended sine kernel.
\par
The determinantal formula due to Okounkov and Reshetikhin (\cite{okounkovreshetikhin}, see also \cite{ferrarispohn}, \cite{borodingorin} or \cite{borodinrains}) is at the center of our proof, although the statement of Theorem \ref{thm1} does not require its knowledge.
\subsection{The determinantal formula}
Okounkov and Reshetikhin have shown in \cite{okounkovreshetikhin} that the pushforward of $\mathbb{P}_q$ onto $\Omega$ is a determinantal point procces. Let us state this fact precisely. We define the kernel $\mathcal{K}_q : E \times E \rightarrow \mathbb{R}$ by :
\begin{align} \label{defKq}
\mathcal{K}_{q}(t_1,h_1;t_2,h_2)= \frac{1}{(2i\pi)^{2}}\int_{|z|=1 \pm \epsilon}\int_{|w|=1 \mp \epsilon} \frac{1}{z-w} \frac{\Phi(t_1,z)}{\Phi(t_2,w)}\frac{dz dw}{z^{h_1+\frac{|t_1|+1}{2}}w^{-h_2-\frac{|t_2|+1}{2}}}
\end{align}
where one picks the plus sign for $t_1 \geq t_2$ and the minus sign otherwise. The function $\Phi$ is defined by :
\begin{align*}
\Phi(t,z)= &\frac{(q^{1/2}/z;q)_\infty}{(q^{1/2+t}z;q)_\infty} \quad \text{for $t \geq 0$} \\
&\frac{(q^{1/2-t}/z;q)_\infty}{(q^{1/2}z;q)_\infty} \quad \text{for $t < 0$},
\end{align*}
where $(x;q)_\infty$ is a $q$ version of the Pochhammer symbol :
\begin{align*}
(x;q)_{\infty}=\prod_{k=0}^{\infty}(1-xq^{k}),
\end{align*}
and $\epsilon$ is a sufficiently small positive number, which allows to avoid the singularities of the ratio :
\begin{align*}
\frac{\Phi(t_1,z)}{\Phi(t_2,w)}.
\end{align*}
Observe that there is no need of defining the square root in formula (\ref{defKq}), since, for $(t,h) \in E$, if there exists a plane partition $\pi$ such that $(t,h) \in \mathfrak{S}(\pi)$, then we have by construction that :
\begin{align*}
h+ \frac{|t|+1}{2} \in \Z.
\end{align*}
Okounkov-Reshetikhin determinantal formula is then the following statement :
\begin{thm}[Okounkov-Reshetikhin, \cite{okounkovreshetikhin}, 2003] \label{thmokresh} For any positive integer $l \in \mathbb{N}$ and any subset $m=\{ (t_1,h_1),...,(t_l,h_l) \} \subset E$, we have :
\begin{align} \label{detformula}
\mathbb{E}_q[c_m(\pi)]=\det \left( \mathcal{K}_q(t_i,h_i;t_j,h_j) \right)_{i,j=1}^l
\end{align}
where $\mathbb{E}_q$ denotes the expectation with respect to $\mathbb{P}_q$.
\end{thm}
The measure $\P_q$ is a particular case of a Schur process. Schur processes, first introduced in \cite{okounkovreshetikhin}, are dynamical generalizations of Schur measures (\cite{okounkovschurmeasures}). For a more elementary treatment of Schur measures and Schur processes, see e.g. \cite{borodinrains}, \cite{borodingorin} and references therein. See also \cite{johansson} for an other approach.
\subsection{The limit process} \label{seclimproc}
In the same article, Okounkov and Reshetikhin proved a scaling limit theorem for $\mathbb{P}_q$, when $r=-\log(q)$ tends to $0^+$, which we now formulate. Let us define :
\begin{align*}A:= \lbrace (t,x) \in \mathbb{R}^{2} , |2\cosh(t/2)-e^{-x}| < 2 \rbrace  \subset \mathbb{R}^2,
\end{align*}
and for $(\tau,\chi) \in A$, let $z(\tau,\chi)$ be the intersection point of the circles $C(0,e^{-\tau/2})$ and $C(1,e^{-\tau/4-\chi/2})$ with positive imaginary part, see figure 1.
\begin{figure}
\includegraphics[scale=0.5,clip=true,trim=0cm 0cm 0cm 0cm]{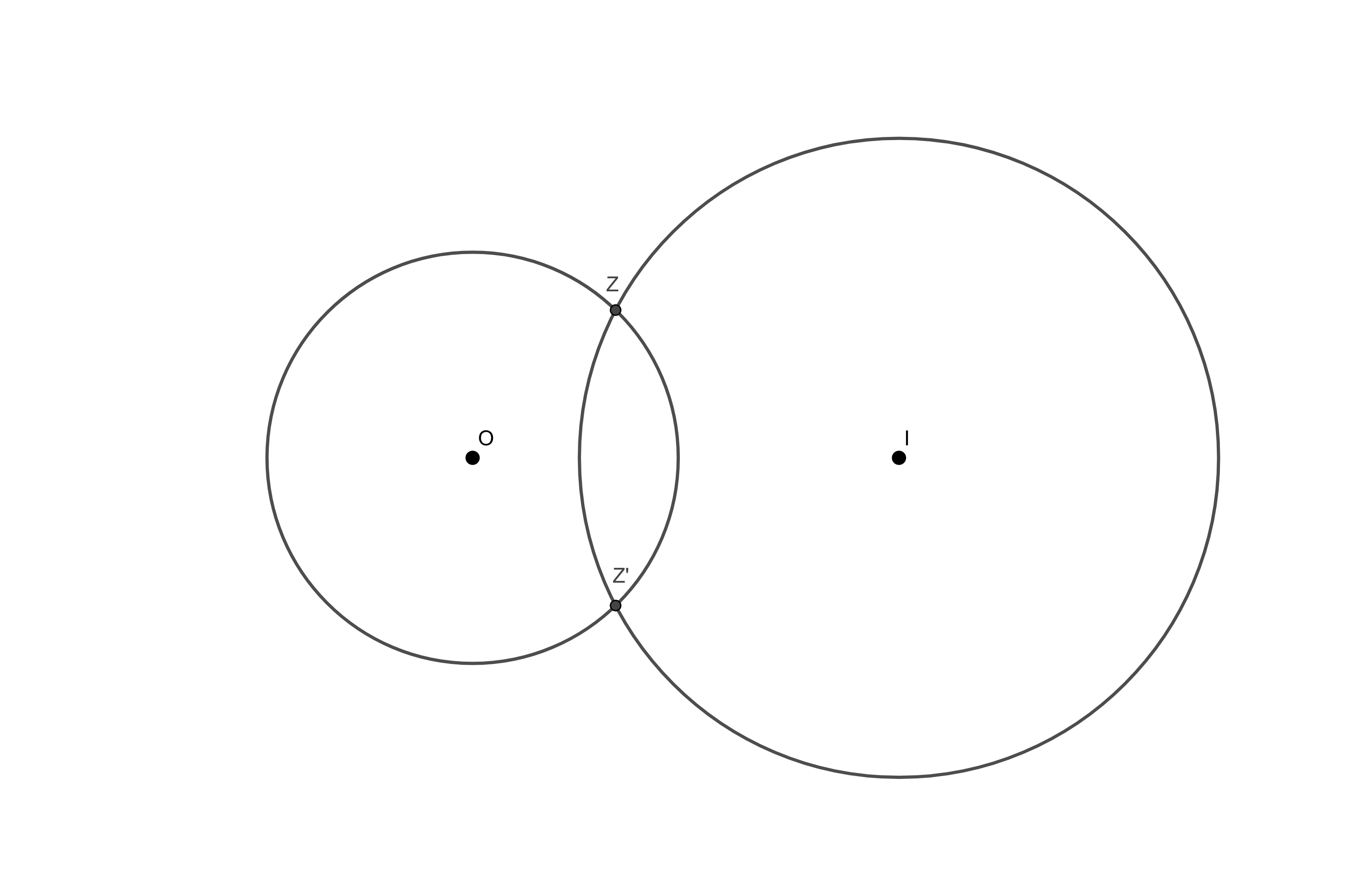}
\caption{The circles $C(0,e^{-\tau/2})$ and $C(1,e^{-\tau/4-\chi/2})$ and their intersection points $z=z(\tau,\chi)$ and its complex conjugate $z'=\overline{z(\tau,\chi)}$.}
\end{figure}
The condition $(\tau,\chi) \in A$ guarantees that $z(\tau,\chi)$ exists and is not real. For $(\tau,\chi) \in A$, we define the translation invariant kernel $\mathcal{S}_{z(\tau,\chi)} : E \rightarrow \mathbb{C}$ :
\begin{align} \label{defdynsin}
\mathcal{S}_{z(\tau,\chi)}(\Delta t, \Delta h)=\frac{1}{2i\pi}\int_{\overline{z(\tau,\chi)}}^{z(\tau,\chi)}(1-w)^{\Delta t}w^{-\Delta h - \frac{\Delta t}{2}} \frac{dw}{w},
\end{align}
where the integration path crosses $(0,1)$ for $\Delta t \geq 0$ and $(-\infty,0)$ for $\Delta t < 0$. For reasons explained below, this kernel will be called the \textit{extended sine kernel}. Then, the following holds :
\begin{thm} [Okounkov-Reshetikhin, \cite{okounkovreshetikhin}, 2003] \label{thmcvokresh} For all $(\tau, \chi) \in A$ and all $m=\{ (t_1,h_1),...,(t_l,h_l) \} \subset E$, we have :
\begin{align*}
\lim_{r\rightarrow 0}  \E_{e^{-r}} \left[c_{\frac{1}{r}(\tau,\chi)+m } \right]= \det \left( \mathcal{S}_{z(\tau,\chi)}(t_i-t_j,h_i-h_j) \right)_{i,j=1}^l.
\end{align*}
\end{thm}
In lemma \ref{lemcv} below, we give the speed of convergence.
\par
For $(\tau, \chi) \in A$, the kernel $\mathcal{S}_{z(\tau,\chi)}$ defines a determinantal point process on $E$, i.e. a probability measure $\P_{(\tau,\chi)}$ on $\Omega$ defined by :
\begin{align} \label{sindynproc}
\forall m=\{ (m^1_1,m^2_1), ..., (m^1_l,m^2_l) \} \subset E, \quad \E_{(\tau,\chi)} [c_m] = \det\left( \mathcal{S}_{z(\tau,\chi)}(m_i^1-m_j^1,m_i^2-m_j^2)\right)_{i,j=1}^l,
\end{align}
where $\E_{(\tau,\chi)}$ is the expectation with respect to $\P_{(\tau,\chi)}$. Theorem \ref{thmcvokresh} means that, when $q$ approaches $1$, if we scale each coordinate of a point process coming from a plane partition distributed according to $\P_q$ by a factor of $r=-\log(q)$ and then zoom around $(\tau,\chi)$, the obtained point process behaves as if it were distributed according to $\P_{(\tau,\chi)}$. This point process can be seen as a two-dimensional or dynamical version of the usual discrete sine-process on $\Z$ (see e.g. \cite{boo} or \cite{borodingorin} and references therein). Indeed, setting $\Delta t =0$ in (\ref{defdynsin}) leads to :
\begin{align*}
\mathcal{S}_{z(\tau,\chi)}(0,\Delta h) = e^{\frac{\tau \Delta h}{2}} \frac{\sin(\phi\Delta h)}{\pi\Delta h}
\end{align*}
where $z(\tau,\chi)= e^{-\frac{\tau}{2}+i\phi}$. Note that we can ignore the factor $e^{\frac{\tau \Delta h}{2}}$ since it will disappear from any determinant of the form (\ref{sindynproc}).
\subsection{Main result}
For $r >0$, we define the set $A_{r} \subset E$ by :
\begin{align*}
A_{r}= r^{-1}A  \cap E.
\end{align*}
For brievety, we write $\P_r$ (resp. $\E_r$) instead of $\P_{e^{-r}}$ (resp. $\E_{e^{-r}}$). For a continuous compactly supported function $f : \R^2 \rightarrow \R$, and a finite subset $m \subset E$, we define the random variable :
\begin{align} \label{sum}
\Sigma(f,m,r)=r^2 \sum_{(t,h) \in A_{r} }f(rt,rh)c_{(t,h)+m}
\end{align}
and the deterministic integral :
\begin{align*}
I(f,m) = \int_{A} f(\tau,\chi)\E_{(\tau,\chi)} [ c_m ] d\tau d\chi.
\end{align*}
Our Theorem establishes that, under $\P_r$, the sum $\Sigma(f,m,r)$ converges to $I(f,m)$. This theorem can thus be interpreted as a weak law of large numbers for functionals of random plane partitions.
\begin{thm} \label{thm1}For every continuous compactly supported function $f: \mathbb{R}^2 \rightarrow \mathbb{R}$, every finite subset $m \subset E$, one has :
\begin{align}
\forall \varepsilon >0, \quad \lim_{r \rightarrow 0} \mathbb{P}_{r} \left( |\Sigma(f,m,r) - I(f,m) | > \varepsilon \right) =0.
\end{align}
\end{thm}
\begin{rem}The assumption of compactness of the support of the function $f$ is used for the uniformity of constants in estimations of averages and variances. It might interesting to see if Theorem \ref{thm1} still holds for a wider class of functions, e.g. Schwartz functions.
\end{rem}
\subsection{Comparison with other models}
In the context of the Plancherel measure on usual partitions, the same theorem was proved by A.I. Bufetov in \cite{bufetovgafa}, lemma 4.4., in order to prove the Vershik-Kerov conjecture concerning the entropy of the Plancherel measure. Indeed, the poissonization of the Plancherel measure is a determinantal point process with the discrete Bessel kernel :
\begin{align*}
\mathcal{J}_\theta(x,y)=\theta \frac{J_x(2\theta)J_{y+1}(2\theta)-J_{x+1}(2\theta)J_y(2\theta)}{x-y}, \quad x,y \in \Z,
\end{align*}
where $J_x(2\theta)$ are the Bessel functions. While the crucial inequality
\begin{align*}
|\mathcal{J}_\theta(x,y)| \leq \frac{C}{|x-y|+1},
\end{align*}
which expresses the decay of correlations for the poissonized Plancherel measure, is here almost immediate, such an inequality in our two-dimensional model requires some efforts. In lemma \ref{lemdecK} below, we prove analogous inequalities for the covariances, which estimate the decay of correlations. We were not able to control the value of the kernel $\mathcal{K}_q$ at different points, but rather for products of the kernel. While the kernel $\mathcal{J}_\theta$ is symmetric, our kernel $\K_q$ is not, and this facts reflects in the need of taking products ; besides, the single value of the kernel of a determinantal process at different points is not always relevant, as one has to consider determinants.
\\
\par
The extended sine kernel appears in many other models as the kernel of the bulk scaling limit of two dimensional statistical mechanics models, for example non-intersecting paths (see e.g. \cite{gorinpaths} and references therein). It also has a continuous counter-part arising in the Dyson's brownian motion model (see e.g. \cite{katoritanemura}). We thus think that a similar law of large numbers should hold for a wide class of discrete determinantal point processes which admit the process with the extended sine-kernel as a limit in the bulk.
\subsection{Organisation of the paper}The paper is organized as follows. In section \ref{sec2}, we introduce notations and state preliminary lemmas, for which the proofs are given later.
\par
Lemma \ref{lemcv} says that the error term in Theorem \ref{thmcvokresh} is of order less than $r$, and as a consequence, the expectation of the sum $\Sigma$ converges to the integral $I$.
\par
Lemma \ref{lemdecK} gives a suitable control on the decay of correlations, and, together with Lemma \ref{lemdiag}, implies that the variance of the sum $\Sigma$ tends to zero.
\par
In section \ref{sec3}, we prove Theorem \ref{thm1}, admitting the lemmas from the preceding section.
\par
Section \ref{sec4} is devoted to the proofs of the lemmas.
\subsection{Acknowledgements} I would like to thank Alexander Bufetov for posing the problem to me and for helpful discussions. I also would like to thank Alexander Boritchev, Nizar Demni and Pascal Hubert for helpful discussions and remarks. This project has received funding from the European Research Council (ERC) under the European Union's Horizon 2020 research and innovation programme (grant agreement N 647133).
\section{Notation and preliminary results} \label{sec2}
\subsection{Notation} For functions  $f,g : [0,+\infty) \rightarrow \C$, we write :
\begin{align*}
f=O(g)
\end{align*}
or :
\begin{align*}
f \lesssim g
\end{align*}
if there exists $C>0$ and $r_0 >0$, such that for all $r \leq r_0$, one has :
\begin{align*}
|f(r)| \leq C|g(r)|.
\end{align*}
We write :
\begin{align*}
f \asymp g
\end{align*}
whenever $f=O(g)$ and $g=O(f)$.\\
\par
The real part of a complex number $z\in \C$ will be denoted by $\mathfrak{R}z$.
\subsection{Preliminary results}
The first lemma we need estimates the error term in Theorem \ref{thmcvokresh}.
\begin{lem} \label{lemcv} For any compact $K \subset \R^2$, and any finite subset $m \subset E$, there exists $C>0$, such that for all $r>0$ sufficiently small and all $(\tau,\chi) \in A \cap K$, one has :
\begin{align}
|\mathbb{E}_r[c_{\frac{1}{r}(\tau,\chi)+m}]-\mathbb{E}_{(\tau,\chi)}[c_m]| \leq Cr.
\end{align}
\end{lem}
The following lemma expresses the decay of correlations of the process $\P_r$.
\begin{lem}\label{lemdecK} Let $K \subset \R^2$ be a compact set, and let $m \subset \E$ be finite. Let $\overline{m}$ denote the supremum norm of $m$ : \begin{align*}
\overline{m}=\max \{ |t|,|h|, \hspace{0.1cm} (t,h) \in m \}.
\end{align*}Then for any $\alpha \in (0,1)$, there exists $C$ which only depends on $\alpha$, $K$ and $m$, such that for all $r >0$ sufficiently small and any $(\tau_1,\chi_1),(\tau_2,\chi_2) \in A\cap K$ such that :
\begin{align} \label{conddist}
\max \{ |\tau_1-\tau_2|, |\chi_1-\chi_2| \} > \overline{m}r,
\end{align}
one has :
\begin{align*}
\left| \E_r \left[ c_{\frac{1}{r}(\tau_1,\chi_1)+m}c_{\frac{1}{r}(\tau_2,\chi_2)+m} \right]-\E_r \left[ c_{\frac{1}{r}(\tau_1,\chi_1)+m}\right]\E_r\left[c_{\frac{1}{r}(\tau_2,\chi_2)+m} \right] \right | \leq \frac{C\exp\left(-r^{-\alpha}\right)}{|\tau_1-\tau_2|^2}
\end{align*}
when $\tau_1 \neq \tau_2$, and :
\begin{align*}
\left| \E_r \left[ c_{\frac{1}{r}(\tau,\chi_1)+m}c_{\frac{1}{r}(\tau,\chi_2)+m} \right]-\E_r \left[ c_{\frac{1}{r}(\tau,\chi_1)+m}\right]\E_r\left[c_{\frac{1}{r}(\tau,\chi_2)+m} \right] \right | \leq \frac{Cr}{|\chi_1-\chi_2|}
\end{align*}
when $\tau_1=\tau_2=\tau$.
\end{lem}
The last lemma we need is obtained as a simple corollary of proposition \ref{propcont} below.
\begin{lem}\label{lemdiag}
For any compact $K \subset \R$ and any finite subsets $m,m' \subset E$, there exists $C$ such that for any $(\tau,\chi) \in A\cap K$, and any sufficiently small $r>0$ :
\begin{align*}
\left| \E_r \left[ c_{\frac{1}{r}(\tau,\chi)+m} c_{\frac{1}{r}(\tau,\chi)+m'}\right] - \E_r \left[ c_{\frac{1}{r}(\tau,\chi)+m}\right] \E_r \left[ c_{\frac{1}{r}(\tau,\chi)+m'} \right] \right| \leq C. 
\end{align*}
\end{lem}
\section{Proof of Theorem \ref{thm1}} \label{sec3}
Let $K \subset \R^2$ be a compact containing the support of $f$. We denote by $A_{r,K}$ the set :
\begin{align*}
A_{r,K}=r^{-1}(A \cap K ) \cap E.
\end{align*}
By lemma \ref{lemcv}, we have :
\begin{align*}
\E_r \Sigma( f,m,r) - r^2 \sum_{(t,h) \in A_{r,K}} f(rt,rh)  \E_{(rt,rh)} \left( c_{m} \right) \lesssim  |A_{r,K}| r^{3},
\end{align*}
where $|A_{r,K}|$ is the cardinality of $A_{r,K}$. We first remark that :
\begin{align} \label{estimArK}
|A_{r,K}| \asymp r^{-2},
\end{align}
which implies :
\begin{align*}
\E_{r} \Sigma (f,m,r)= r^2 \sum_{(t,h) \in A_{r,K}} f(rt,rh) \E_{(rt,rh)} \left(c_{m}\right) + O(r).
\end{align*}
Observing then that 
\begin{align*}
r^2 \sum_{(t,h) \in A_{r,K}} f(rt,rh) \E_{(rt,rh)} \left(c_{m}\right)
\end{align*}
is a Riemann sum for the integral $I(f,m)$, we obtain that :
\begin{align*}
\lim_{r \rightarrow 0} \E_{r} \Sigma(f,m,r) = I(f,m).
\end{align*}
By the Chebyshev inequality, it suffices now to prove that :
\begin{align} \label{var0}
\text{Var}_{r} \left( \Sigma(f,m,r) \right) \rightarrow 0 
\end{align}
as $r$ tends to $0$, where :
\begin{multline} \label{variance}
\text{Var}_{r} \left( \Sigma(f,m,r) \right) = \E_r \left[ \left(\Sigma(f,m,r) - \E_r \Sigma(f,m,r) \right)^2 \right] \\
=r^4 \sum_{(t_1h_1),(t_2,h_2) \in A_{r,K}}  f(rt_1,rh_1)f(rt_2,rh_2) \\
\times \left(\E_r ( c_{(t_1,h_1)+m}c_{(t_2,h_2)+m} ) - \E_r ( c_{(t_1,h_1)+m} ) \E_r ( c_{(t_2,h_2)+m} ) \right).
\end{multline}
We set $ \overline{m}=\max\{ |t|,|h|, \hspace{0.1cm} (t,h) \in m \}$, and we partition $A_{r,K}^2$ into three sets :
\begin{align*}
A_{r,K}^2=A_{r,k}^{>} \sqcup A_{r,K}^{>=} \sqcup A_{r,K}^{\leq},
\end{align*}
where :
\begin{align*}
A_{r,k}^{>} &=\left\{ (t_1,h_1), (t_2,h_2) \in A_{r,K}, \hspace{0.1cm} \max\{ |t_1-t_2|, |h_1-h_2|\} > \overline{m}, \hspace{0.1cm} t_1 \neq t_2 \right\}, \\
A_{r,K}^{>=} &= \left\{ (t,h_1), (t,h_2) \in A_{r,K}, \hspace{0.1cm} |h_1-h_2| > \overline{m} \right\},\\
A_{r,K}^{\leq} &= A_{r,K}^2 \setminus \left(A_{r,k}^{>} \sqcup A_{r,K}^{>=} \right) = \left\{ (t_1,h_1),(t_2,h_2) \in A_{r,K}, \hspace{0.1cm} \max\{ |t_1-t_2|, |h_1-h_2|\} \leq \overline{m} \right\}.
\end{align*}
We first estimate the variance (\ref{variance}) by :
\begin{multline} \label{estimvar}
\text{Var}_{r} \left( \Sigma(f,m,r) \right) 
\leq C r^4 \left( \sum_{\left((t_1h_1),(t_2,h_2)\right) \in A_{r,K}^>}  \left|\E_r ( c_{(t_1,h_1)+m}c_{(t_2,h_2)+m} ) - \E_r ( c_{(t_1,h_1)+m} ) \E_r ( c_{(t_2,h_2)+m} ) \right| \right. \\
+ \sum_{ \left((t,h_1),(t,h_2)\right) \in A_{r,K}^{>=} } \left|\E_r ( c_{(t,h_1)+m}c_{(t,h_2)+m} ) - \E_r ( c_{(t,h_1)+m} ) \E_r ( c_{(t,h_2)+m} ) \right| \\
\left. + \sum_{\left((t_1h_1),(t_2,h_2)\right) \in A_{r,K}^{\leq}}  \left|\E_r ( c_{(t_1,h_1)+m}c_{(t_2,h_2)+m} ) - \E_r ( c_{(t_1,h_1)+m} ) \E_r ( c_{(t_2,h_2)+m} ) \right|\right),
\end{multline}
where $C$ only depends on $f$. Let $(t_1,h_1),(t_2,h_2) \in A_{r,K}$. By definition, there exists $(\tau_1,\chi_1),(\tau_2,\chi_2) \in A \cap K \cap rE$ such that :
\begin{align*}
(t_1,h_1)=\frac{1}{r}(\tau_1,\chi_1), \hspace{0.1cm} (t_2,h_2) =\frac{1}{r}(\tau_2,\chi_2).
\end{align*}
We first consider the case when $\left((t_1,h_,),(t_2,h_2)\right) \in A_{r,K}^>$. The corresponding points $(\tau_1,\chi_1),(\tau_2,\chi_2)$ satisfy condition (\ref{conddist}), and thus by lemma \ref{lemdecK}, we have in particular the estimate :
\begin{align*}
\left| \E_r \left[ c_{\frac{1}{r}(\tau_1,\chi_1)+m}c_{\frac{1}{r}(\tau_2,\chi_2)+m} \right]-\E_r \left[ c_{\frac{1}{r}(\tau_1,\chi_1)+m}\right]\E_r\left[c_{\frac{1}{r}(\tau_2,\chi_2)+m} \right] \right | \leq Cr,
\end{align*}
where $C$ is uniform. Since :
\begin{align*}
A_{r,K}^> \asymp r^{-4},
\end{align*}
we obtain that :
\begin{align} \label{estimcov1}
\sum_{\left((t_1h_1),(t_2,h_2)\right) \in A_{r,K}^>}  \left|\E_r ( c_{(t_1,h_1)+m}c_{(t_2,h_2)+m} ) - \E_r ( c_{(t_1,h_1)+m} ) \E_r ( c_{(t_2,h_2)+m} ) \right| \leq Cr^{-3}
\end{align}
where $C$ only depends on $f$ and $m$. 
\par
In the case when $\left((t_1,h_,),(t_2,h_2)\right) \in A_{r,K}^{>=}$, we have by lemma \ref{lemdecK} that :
\begin{align*}
\left| \E_r \left[ c_{\frac{1}{r}(\tau_1,\chi_1)+m}c_{\frac{1}{r}(\tau_2,\chi_2)+m} \right]-\E_r \left[ c_{\frac{1}{r}(\tau_1,\chi_1)+m}\right]\E_r\left[c_{\frac{1}{r}(\tau_2,\chi_2)+m} \right] \right | \leq C
\end{align*}
where $C$ is uniform. Since :
\begin{align*}
|A_{r,K}^{<=}| \asymp r^{-3}, 
\end{align*}
we have :
\begin{align}\label{estimcov2}
\sum_{(t,h_1),(t,h_2) \in A_{r,K}^{>=}}  \left|\E_r ( c_{(t,h_1)+m}c_{(t,h_2)+m} ) - \E_r ( c_{(t,h_1)+m} ) \E_r ( c_{(t,h_2)+m} ) \right| \leq Cr^{-3},
\end{align}
where $C$ only depends on $K$ and $m$.
\par 
When $\left((t_1,h_1),(t_2;h_2)\right) \in A_{r,K}^{\leq}$, there exists finite subsets $m',m'' \subset E$ and $(\tau,\chi) \in A\cap K$ such that :
\begin{align*}
(t_1,h_1)+m= \frac{1}{r}(\tau,\chi) +m', \hspace{0.1cm} (t_2,h_2)+m = \frac{1}{r}(\tau,\chi) +m''.
\end{align*}
Observe that there are only a finite number of possible sets $m'$ and $m''$.Thus, by lemma \ref{lemdiag}, we have :
\begin{align*}
\left| \E_r \left[ c_{\frac{1}{r}(\tau,\chi)+m'} c_{\frac{1}{r}(\tau,\chi)+m''}\right] - \E_r \left[ c_{\frac{1}{r}(\tau,\chi)+m'}\right] \E_r \left[ c_{\frac{1}{r}(\tau,\chi)+m''} \right] \right| \leq C
\end{align*}
where $C$ is uniform. Since :
\begin{align*}
|A_{r,K}^{\leq} |\asymp r^{-2},
\end{align*}
we have :
\begin{align} \label{estimcov3}
\sum_{(t,h_1),(t,h_2) \in A_{r,K}^{\leq} }\left|   \E_r ( c_{(t,h_1)+m}c_{(t,h_2)+m} ) - \E_r ( c_{(t,h_1)+m} ) \E_r ( c_{(t,h_2)+m} ) \right|  \leq Cr^{-2}
\end{align}
Thus, recalling the estimation (\ref{estimvar}), the inequalities (\ref{estimcov1}), (\ref{estimcov2}) and (\ref{estimcov3}) establish (\ref{var0}). Theorem \ref{thm1} is proved, assuming lemmas \ref{lemcv}, \ref{lemdecK} and \ref{lemdiag}.
\section{Proof of lemmas \ref{lemcv}, \ref{lemdecK} and  \ref{lemdiag}} \label{sec4}
\subsection{Proof of lemma \ref{lemcv}}We here follow the proof of \cite{okounkovreshetikhin}, giving the error terms in the asymptotics we use. We define the dilogarithm function as being the analytic continuation of the series :
\begin{align*}
\text{dilog}(1-z)=\sum_{n \geq 1} \frac{z^n}{n^2}, \quad |z|<1,
\end{align*}
with a cut along the half-line $(1,+\infty)$. We first prove that :
\begin{align} \label{estimdilog}
-\log(z,q)_\infty = r^{-1}\text{dilog}(1-z) + O(1).
\end{align}
as $q=e^{-r}$ tends to $1^-$. Indeed, we have :
\begin{align*}
\log(z,q)_{\infty}&=\sum_{k \geq 0} \log(1-zq^k) = -\sum_{k \geq 0} \sum_{n \geq 1} \frac{z^nq^{nk}}{n} \\
&=-\sum_{n \geq 1} \frac{z^n}{n} \sum_{k \geq 0} q^{nk} =-\sum_{n \geq 0} \frac{z^n}{n} \frac{1}{1-q^n}.
\end{align*}
With $q=e^{-r}$, we have :
\begin{align*}
\frac{r}{1-e^{-nr}}=\frac{1}{n -n^2r + ...}=\frac{1}{n}(1+nr+...)
\end{align*}
and thus :
\begin{align*}
\left| \frac{z^n}{n}\left(\frac{1}{1-e^{-nr}}- \frac{1}{n}\right)\right| \leq rn|z|^n,
\end{align*}
which establishes (\ref{estimdilog}). 
\par
Let $K \subset \R^2$ be compact and let $(\tau,\chi) \in A\cap K$. We assume that $\tau \geq 0$, see \ref{secrk} below for the case $\tau \leq 0$.  We introduce the function :
\begin{align*}
S(z;\tau,\chi)=-(\tau/2+ \chi)\log(z)-\dilog(1-1/z)+ \dilog(1-e^{-\tau}z),
\end{align*}
and denote by  $\gamma_\tau$ the circle :
\begin{align*}
\gamma_\tau = \{ z \in \C, \hspace{0.1cm} |z|=e^{\tau/2} \}.
\end{align*}
By the estimate (\ref{estimdilog}) and formula (\ref{defKq}), we have that, for all $z$ and $w$ sufficiently closed to $\gamma_\tau$ : 
\begin{multline} \displaystyle
\left| \frac{\Phi(\tau/r+ t_1,z)}{\Phi(\tau/r+t_2,w)}\frac{1}{z^{\chi/r+h_1+\tau/2r+(t_1+1)/2}w^{-\chi/r-h_2-\tau/2r-(t_2+1)/2}} \right| \\
=  \exp \left( \frac{1}{r} \left( \mathfrak{R} S(z;\tau,\chi)-\mathfrak{R}S(w;\tau,\chi) \right) + O(1)  \right) ,
\end{multline}
where the $O(1)$ term only depens on $K$, $(t_1,h_1)$ and $(t_2,h_2)$. An observation made in \cite{okounkovreshetikhin} states that the real part of $S$ on the circle $\gamma_\tau$ is constant, namely :
\begin{align} \label{reS}
\mathfrak{R}S(z;\tau,\chi)=-\frac{\tau}{2}(\tau/2+\chi),\quad z \in \gamma_\tau.
\end{align}
It is also shown in \cite{okounkovreshetikhin} that, since $(\tau,\chi) \in A$, the function $S$ has two distinct critical points on $\gamma_\tau$ : $e^{\tau}z(\tau,\chi)$ and its complex conjugate. The computation of the gradient of the real part of $S$ on $\gamma_\tau$ lead then the authors of \cite{okounkovreshetikhin} to deform the circle $\gamma_\tau$ into simple contours $\gamma_\tau^>$ and $\gamma_\tau^<$, both crossing the two critical points and verifying :
\begin{equation}
\begin{split}
z \in \gamma_\tau^> \Rightarrow \mathfrak{R}S(z;\tau,\chi) \geq -\frac{\tau}{2}(\tau/2+\chi), \\
z \in \gamma_\tau^< \Rightarrow \mathfrak{R}S(z;\tau,\chi) \leq -\frac{\tau}{2}(\tau/2+\chi),
\end{split}
\end{equation}
with equality only for $z\in \left\{e^\tau z(\tau,\chi), e^\tau \overline{z(\tau,\chi)} \right\}$, see figure 2.
\begin{figure}[!h] \label{fig1}
\centering
\includegraphics[scale=0.6,clip=true,trim=0cm 0cm 0cm 0cm]{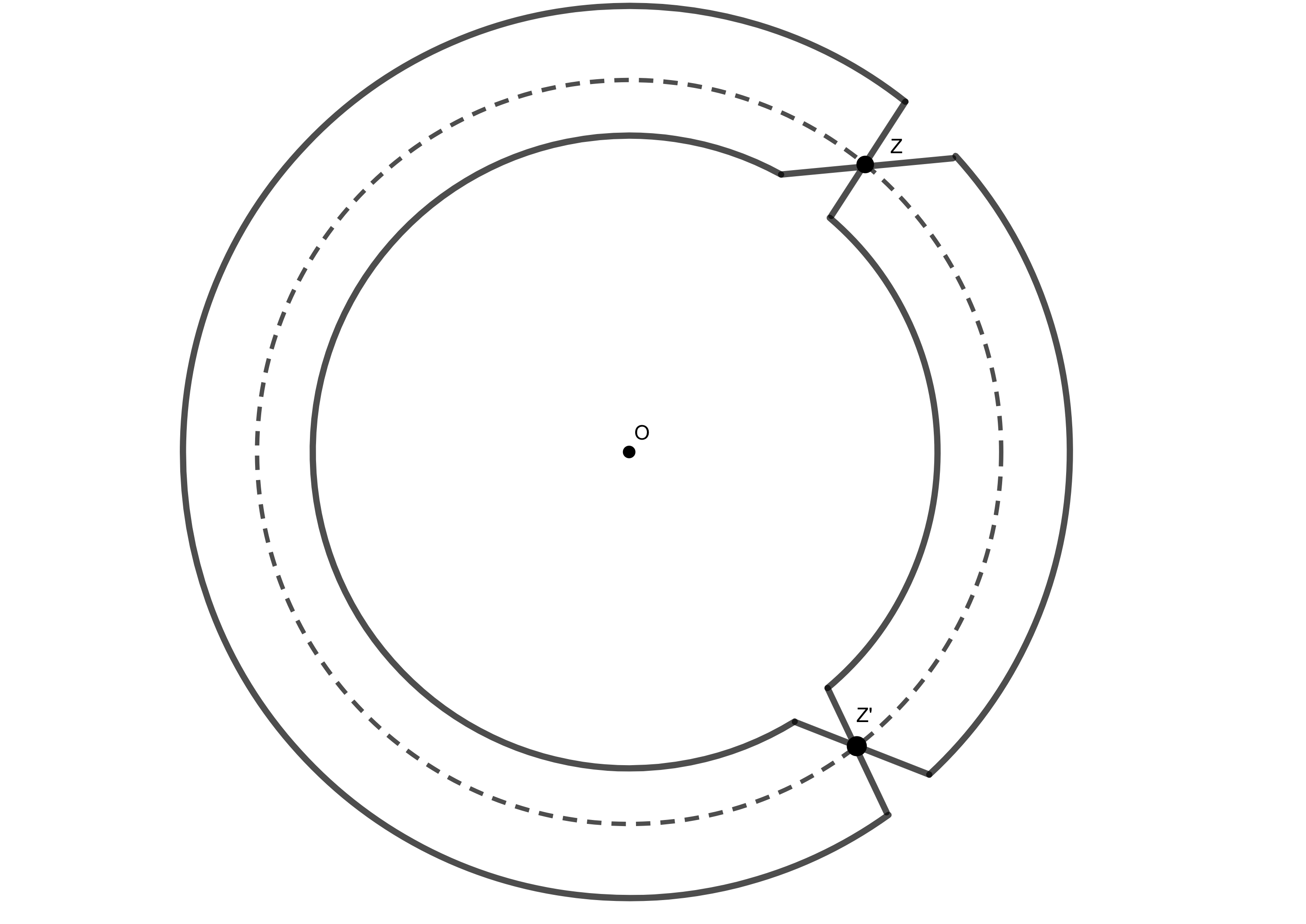}
\caption{The contours $\gamma_\tau^>$ and $\gamma_\tau^<$ are the thick contours and the circle $\gamma_\tau$ is the dotted circle.}
\end{figure}
These simple facts imply that the integral :
\begin{align} \label{intdef}
 \int_{z \in \gamma_\tau^<} \int_{w \in \gamma_\tau^>} \exp \left( \frac{1}{r} \left( S(z;\tau,\chi)-S(w;\tau,\chi) \right)\right) \frac{dzdw}{z-w}
\end{align}
goes to zero as $r$ tends to zero. Actually, the dominated convergence theorem implies that the integral (\ref{intdef}) is $O\left(\exp\left(-r^{-\alpha}\right)\right)$ for any $\alpha \in (0,1)$.
\par 
Picking the residue at $z=w$, we arrive at :
\begin{multline} \label{sumKq}
\mathcal{K}_{e^{-r}} \left( \frac{\tau}{r} + t_1, \frac{\chi}{r} +h_1 ; \frac{\tau}{r} + t_2, \frac{\chi}{r} +h_2\right) = \frac{1}{(2i\pi)^2} \int_{z \in \gamma_\tau^<} \int_{w \in \gamma_\tau^>} \exp \left( \frac{1}{r} \left( S(z;\tau,\chi)-S(w;\tau,\chi) \right)+ O(1) \right) \frac{dzdw}{z-w} \\
+ \frac{1}{2i\pi} \int_{e^\tau \overline{z(\tau,\chi)}}^{e^\tau z(\tau,\chi)} \frac{(q^{1/2+\tau/r+t_2}w;q)_\infty}{(q^{1/2+ \tau/r +t_1};q)_\infty} \frac{dw}{w^{h_1-h_2+(t_1-t_2)/2}},
\end{multline}
where the path of integration for the second integral crosses the interval $(0,e^\tau)$ for $t_1\geq t_2$ and the half-line $(-\infty,0)$ otherwise. By the preceding discussion, the first integral rapidly tends to zero. Observe now that :
\begin{align*}
\frac{(q^{1/2+\tau/r+t_2}w;q)_\infty}{(q^{1/2+ \tau/r +t_1};q)_\infty} = (1+O(r))\left(1-e^{-\tau}w\right)^{t_1-t_2}
\end{align*}
where the $O(r)$ term only depends on $K$, $t_1$ and $t_2$. Performing the change of variable $w \mapsto e^{-\tau}w$ in the second integral of (\ref{sumKq}), we arrive at :
\begin{align*} \displaystyle
\frac{1}{2i\pi} \int_{e^\tau \overline{z(\tau,\chi)}}^{e^\tau z(\tau,\chi)} \frac{(q^{1/2+\tau/r+t_2}w;q)_\infty}{(q^{1/2+ \tau/r +t_1};q)_\infty} \frac{dw}{w^{h_1-h_2+(t_1-t_2)/2}} = \left( 1 +O(r) \right) e^{-\tau(h_1-h_2-(t_1-t_2)/2)} \mathcal{S}_{\tau,\chi }( t_1-t_2,h_1-h_2).
\end{align*}
The factor $e^{-\tau(h_1-h_2-(t_1-t_2)/2)}$ can be ignored, since it disappears from any determinant of the form (\ref{defdynsin}). Lemma \ref{lemcv} is proved. $\square$.
\subsection{A remark and a proposition} \label{secrk}
For $\tau <0$, one has to replace the function $S$ by :
\begin{align*}
\tilde{S}(z,\tau,\chi)=-(|\tau|/2+\chi)\log(z)-\dilog(1-z)+\dilog(1-e^{-|\tau|}/z).
\end{align*}
The function $\tilde{S}$ innerhits the same properties than the function $S$ : it is constant on the circle $\gamma_{|\tau|}$ and has two complex conjugated critical points on this circle provided $(\tau,\chi) \in A$. This is why we will only consider positve values of $\tau$ in the sequel.\\
\par
The critical points of $S$ are the roots of the quadratic polynomial :
\begin{align*}
(1-1/z)(1-e^{-\tau}z)=e^{-\tau/2-\chi}.
\end{align*}
For this reason, we have the following proposition :
\begin{prop} \label{propcont}
For any fixed $(\Delta t, \Delta h) \in E$, the function :
\begin{align*}
(\tau,\chi) \mapsto \mathcal{S}_{\tau,\chi}(\Delta t, \Delta h)
\end{align*}
is continuous.
\end{prop}
\subsection{Proof of lemma \ref{lemdecK}}
Let $m \subset E$ be finite, of cardinality $l$, let $K \subset \R^2$ be a compact set and let $(\tau_1,\chi_1), (\tau_2,\chi_2) \in A\cap K$ be as in the statement of the lemma. The condition (\ref{conddist}) implies that the sets $\frac{1}{r}(\tau_1,\chi_1)+m$ and $\frac{1}{r}(\tau_2,\chi_2)+ m$ are disjoints. Thus, the expectation :
\begin{align*}
\E_r \left[ c_{\frac{1}{r}(\tau_1,\chi_1)+m} c_{\frac{1}{r}(\tau_2,\chi_2)+m} \right]
\end{align*}
is a determinant of size $2l$. In the expansion of the determinant as an alternate sum over permutations of size $2l$, one considers the permutations leaving the sets $\frac{1}{r}(\tau_1,\chi_1)+m$ and $\frac{1}{r}(\tau_2,\chi_2)+m$ invariant. The alternate sum over all such permutations is nothing but the product of determinants :
\begin{align*}
\E_r \left[ c_{\frac{1}{r}(\tau_1,\chi_1)+m} \right ] \E_r \left[ c_{\frac{1}{r}(\tau_2,\chi_2)+m} \right].
\end{align*}
Thus, the terms of the remaining sum :
\begin{align*}
\E_r \left[ c_{\frac{1}{r}(\tau_1,\chi_1)+m} c_{\frac{1}{r}(\tau_2,\chi_2)+m} \right] - \E_r \left[ c_{\frac{1}{r}(\tau_1,\chi_1)+m} \right ] \E_r \left[ c_{\frac{1}{r}(\tau_2,\chi_2)+m} \right]
\end{align*}
all involve factors of the form :
\begin{align} \label{prodK}
\K_{e^{-r}}\left(\frac{1}{r}(\tau_1,\chi_1)+m_{i_1}; \frac{1}{r}(\tau_2,\chi_2)+m_{j_1}\right)\K_{e^{-r}}\left(\frac{1}{r}(\tau_2,\chi_2)+m_{j_2}; \frac{1}{r}(\tau_1,\chi_1)+m_{i_2}\right).
\end{align}
By lemma \ref{lemcv} and proposition \ref{propcont}, the other factors are bounded by a bound only depending on $K$ and $m$. By similar methods as in the proof of Lemma \ref{lemcv}, we will show that the product (\ref{prodK}) is small. 
\par 
The product (\ref{prodK}) can be written as a quadruple integral :
\begin{multline*} \displaystyle
\K_{e^{-r}}\left(\frac{1}{r}(\tau_1,\chi_1)+m_{i_1}; \frac{1}{r}(\tau_2,\chi_2)+m_{j_1}\right)\K_{e^{-r}}\left(\frac{1}{r}(\tau_2,\chi_2)+m_{j_2}; \frac{1}{r}(\tau_1,\chi_1)+m_{i_2}\right) \\
= \frac{1}{(2i\pi)^4} \int_{z \in (1+\varepsilon)\gamma_{\tau_1}} dz\int_{w \in (1-\varepsilon)\gamma_{\tau_2}} dw\int_{z' \in (1-\varepsilon)\gamma_{\tau_2}}dz' \int_{w' \in (1+\varepsilon)\gamma_{\tau_1}} dw' \\
\frac{\exp\left(\frac{1}{r}\left(S(z;\tau_1,\chi_1)-S(w;\tau_2,\chi_2)+S(z';\tau_2,\chi_2)-S(w';\tau_1,\chi_1)\right)+O(1)\right)}{(z-w)(z'-w')}
\end{multline*}
We first consider the case when $\tau_1 \neq \tau_2$, and by symmetry, we assume that $\tau_1 > \tau_2$. One can then deform the contours as previously. Precisely, we now integrate over :
\begin{align*}
z \in \gamma_{\tau_1}^<, \hspace{0.1cm} w \in \gamma_{\tau_2}^>, \hspace{0.1cm} z' \in \gamma_{\tau_2}^<, \hspace{0.1cm} w' \in \gamma_{\tau_1}^>
\end{align*}
in order to have
\begin{align*}
\mathfrak{R}\left(S(z;\tau_1,\chi_1)-S(w';\tau_1,\chi_1) \right) <0, \quad \text{and} \quad \mathfrak{R}\left(S(z';\tau_2,\chi_2)-S(w;\tau_2,\chi_2) \right) <0,
\end{align*}
see figure 3.
\begin{figure}[!h] \label{fig2}
\centering
\includegraphics[scale=0.3,clip=true,trim=0cm 0cm 0cm 0cm]{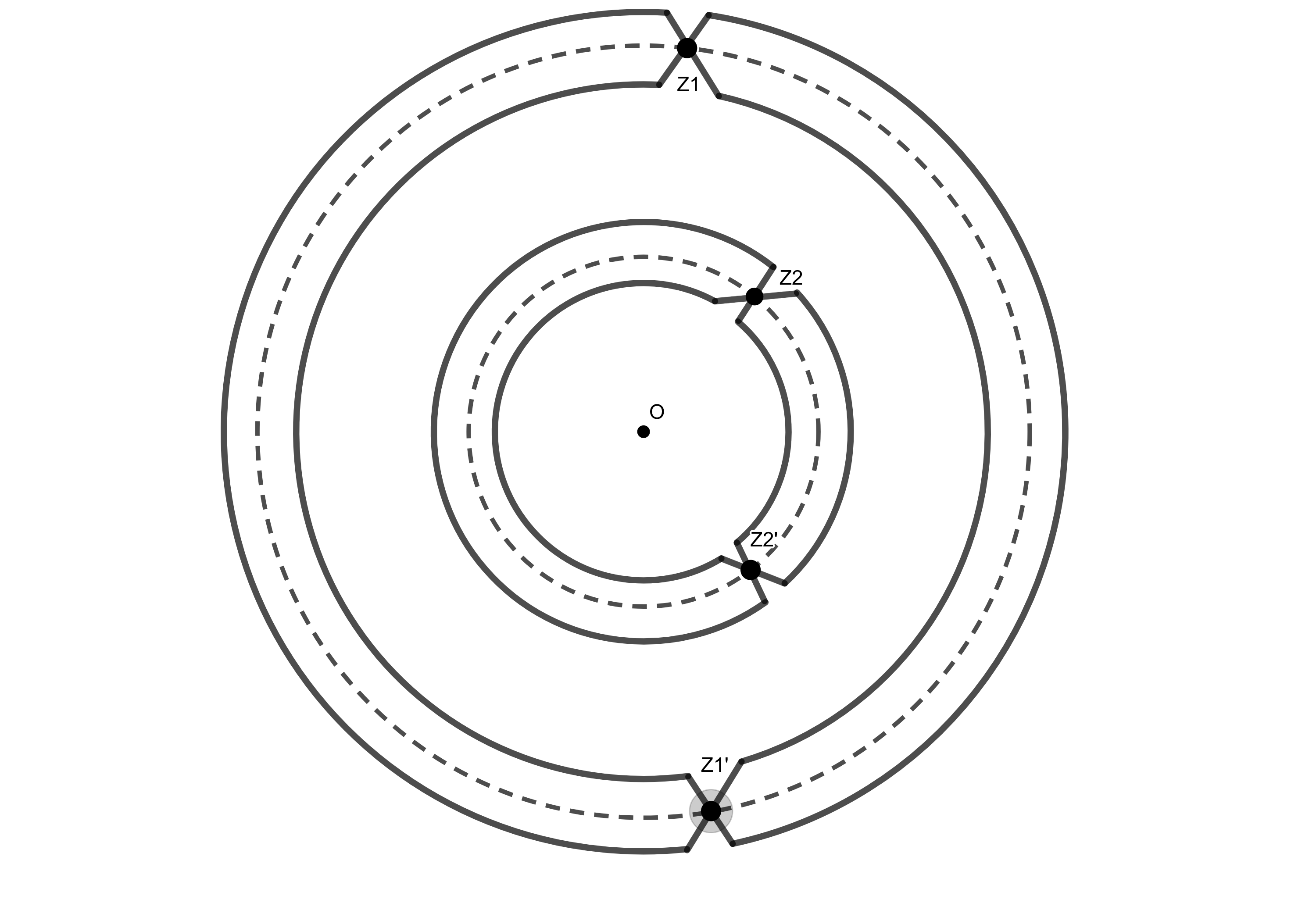}
\caption{The contours $\gamma_{\tau_1}^>$, $\gamma_{\tau_1}^<$ are the thick contours near the dotted circle with the largest radius, the circle $\gamma_{\tau_1}$ ; the contours $\gamma_{\tau_2}^>$ and $\gamma_{\tau_2}^<$ are the thick contours near the dotted circle with the smallest radius $\gamma_{\tau_2}$.}
\end{figure}
These deformations do not affect the value of the integrals, because they involve separate variables. Since, for all $\alpha \in (0,1)$, we have :
\begin{align*}
\frac{\exp\left(\frac{1}{r}\left(S(z;\tau_1,\chi_1)-S(w;\tau_2,\chi_2)+S(z';\tau_2,\chi_2)-S(w';\tau_1,\chi_1)\right)\right)}{\exp(r^{-\alpha})} \rightarrow 0
\end{align*}
as $r \rightarrow 0$, for all $z,z',w,w'$ in the new contours except at a finite number of points, and since 
\begin{align*}
\left| \frac{\exp\left(\frac{1}{r}\left(S(z;\tau_1,\chi_1)-S(w;\tau_2,\chi_2)+S(z';\tau_2,\chi_2)-S(w';\tau_1,\chi_1)\right)\right)}{(z-w)(z'-w')} \right| \leq \frac{1}{|e^{\tau_1}-e^{\tau_2}|^2}
\end{align*}
for all $z,z',w,w'$ in the new contours except at a finite number of points, this concludes the proof for $\tau_1 >\tau_2$.
\par 
For $\tau_1=\tau_2=\tau$, the proof is as follows. One deform the contours as for the preceding case, but now, the deformations affect the value of the kernel since we can avoid the residues at $z=w$ and $z'=w'$. We have for example the following case :
\begin{equation} \label{prodKq2}
\begin{split}
\K&_{e^{-r}}\left(\frac{1}{r}(\tau,\chi_1)+(t_1^1,h_1^1); \frac{1}{r}(\tau,\chi_2)+(t_2^1,h_2^1)_1\right)\K_{e^{-r}}\left(\frac{1}{r}(\tau,\chi_2)+(t_1^2,h_1^2)_2; \frac{1}{r}(\tau,\chi_1)+(t_2^2,h_2^2)_2\right) \\
&=(1+O(1))\left( \frac{1}{(2i\pi)^2}\int_{z \in \gamma_\tau^{<,1}}\int_{w \in \gamma_\tau^{>,2}}... +\frac{1}{2i\pi}\int_w Res_{z=w} f(z,w;\tau,\chi_1,\chi_2)dw \right) \\
&\times
\left(\frac{1}{(2i\pi)^2} \int_{z' \in \gamma_\tau^{<,2}}\int_{w' \in \gamma_\tau^{>,1}}... + \frac{1}{2i\pi}\int_{w'} Res_{z'=w'}f(z',w';\tau,\chi_2,\chi_1)dw' \right),
\end{split}
\end{equation} 
where :
\begin{equation} \label{res}
\begin{split}
\int_w Res_{z=w}f(z,w;\tau,\chi_1,\chi_2)dw=\int_{|w|=e^{\tau/2}, \hspace{0.1cm} |\arg(w)|<\phi_{\tau,\chi_2}} \frac{(q^{1/2+\tau/r+t_2}w;q)_\infty}{(q^{1/2+\tau/r+t_1}w;q)_\infty}\frac{dw}{w^{1/r(\chi_1-\chi_2)+h_1^1-h_2^1 + t_1^1-t_2^1}}, \\
\int_{w'} Res_{z'=w'}f(z',w';\tau,\chi_2,\chi_1)dw=\int_{|w'|=e^{\tau/2}, \hspace{0.1cm} |\arg(w')|<\phi_{\tau,\chi_2}} \frac{(q^{1/2+\tau/r+t_2}w';q)_\infty}{(q^{1/2+\tau/r+t_1}w';q)_\infty}\frac{dw'}{w'^{1/r(\chi_2-\chi_1)+ h_1^2-h_2^2 +  t_1^2-t_2^2}},
\end{split}
\end{equation}
the argument $\phi_{\tau,\chi_2}$ being an argument of $z(\tau,\chi_2)$, see figure 4. Equality (\ref{prodKq2}) is valid when :
\begin{align*}
t_1^1 \geq t_2^1, \hspace{0.1cm} t_1^2 \geq t_2^2, \hspace{0.1cm} \text{and } \arg \left( z(\tau,\chi_1) \right) <\arg \left( z(\tau,\chi_2) \right),
\end{align*}
and the other cases can be treated in a similar way.
\begin{figure}[!h]
\centering
\includegraphics[scale=0.5,clip=true,trim=0cm 0cm 0cm 0cm]{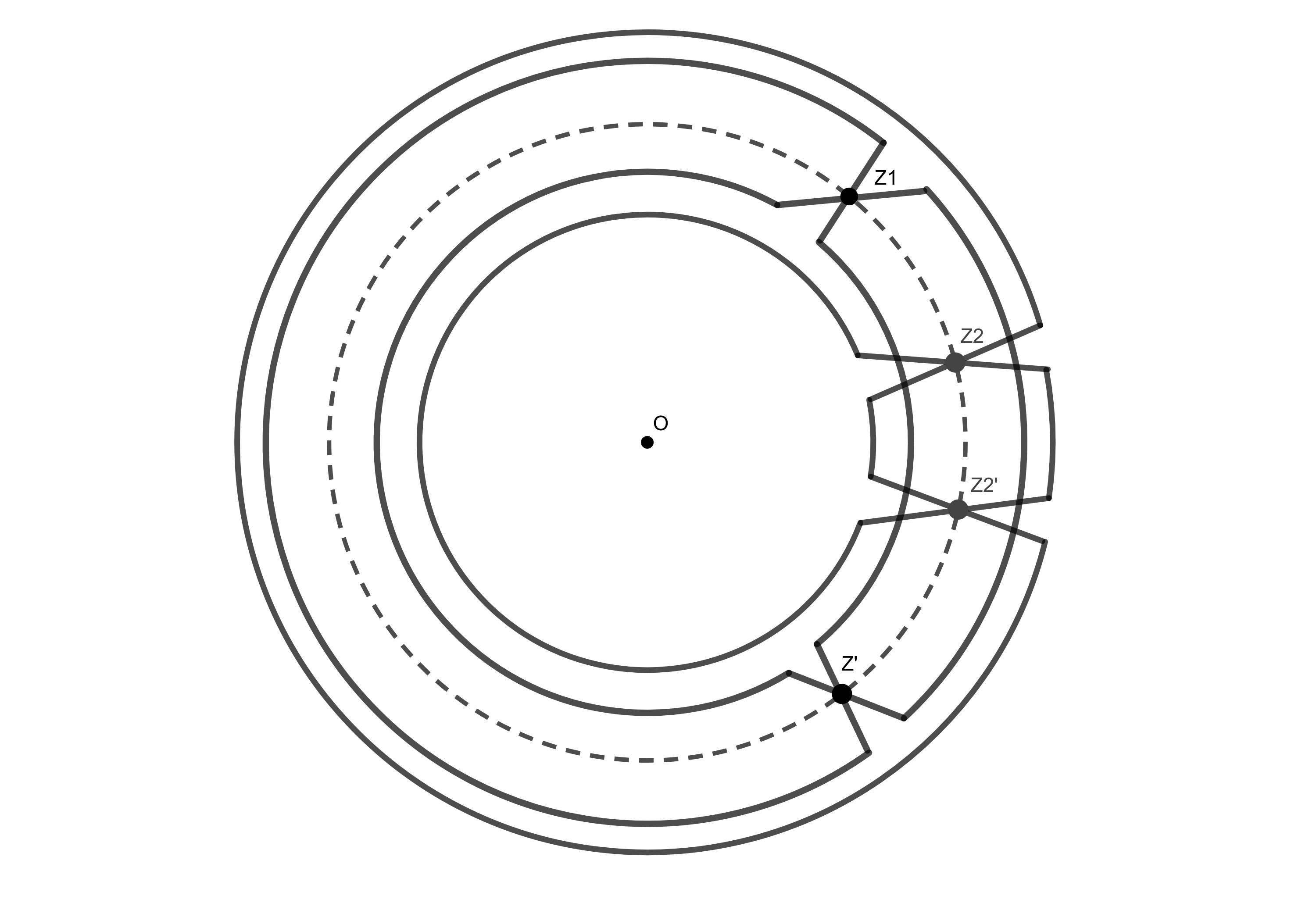}
\caption{The thick contours $\gamma_{\tau}^{>,1}$ and $\gamma_\tau^{<,1}$ cross the dotted circle $\gamma_\tau$ at points $z1=e^\tau z(\tau,\chi_1)$ and $z1'=e^\tau\overline{z(\tau,\chi_1)}$, while the thick contours $\gamma_{\tau}^{>,2}$ and $\gamma_{\tau}^{<,2}$ cross the circle $\gamma_\tau$ at points $z2=e^\tau z(\tau,\chi_2)$ and $z2'=e^\tau \overline{z(\tau,\chi_2)}$.}
\end{figure}
Note that the factor :
\begin{align*}
f_r(w):=\frac{(q^{1/2+\tau/r+t_2}w;q)_\infty}{(q^{1/2+\tau/r+t_1}w;q)_\infty}
\end{align*}
is bounded, since it tends to :
\begin{align*}
(1-e^{-\tau}w)^{t_1-t_2}
\end{align*}
as $r$ tends to $0$. Integrating (\ref{res}) by parts leads to :
\begin{equation} \label{ineqres}
\begin{split}
\left| \int_w Res_{z=w}f(z,w;\tau,\chi_1,\chi_2)dw \right| &\leq  C \frac{\exp\left(\frac{\tau}{2r}(\chi_2-\chi_1)\right)}{1/r|\chi_1-\chi_2|}\\
&+\frac{\exp\left(\frac{\tau}{2r}(\chi_2-\chi_1)\right)}{1/r|\chi_1-\chi_2|}\left| \int_{|w|=1, \hspace{0.1cm} \arg(w)<\phi_{\tau,\chi_2}}f_r'(e^{-\tau/2}w)\frac{dw}{w^{1/r(\chi_1-\chi_2)+\Delta h + \Delta t-1}}\right| \\
&  \leq C\frac{r}{|\chi_1-\chi_2|}\exp \left( \frac{\tau}{2r}(\chi_2-\chi_1)\right),
\end{split}
\end{equation}
and :
\begin{equation} \label{ineq2}
\left| \int_{w'} Res_{z'=w'}f(z',w';\tau,\chi_2,\chi_1)dw \right|  \leq C\frac{r}{|\chi_1-\chi_2|}\exp \left( \frac{\tau}{2r}(\chi_1-\chi_2)\right).
\end{equation}
It is clear that, by construction, we have :
\begin{equation}\label{ineqint}
\left|\int_{z \in \gamma_\tau^{<,1}}\int_{w \in \gamma_\tau^{>,2}}... \right| \leq C \exp \left( \frac{\tau}{2r}(\chi_2-\chi_1)\right),
\end{equation}
and :
\begin{equation}\label{ineqint2}
\left|\int_{z' \in \gamma_\tau^{<,2}}\int_{w' \in \gamma_\tau^{>,1}}... \right| \leq C \exp \left( \frac{\tau}{2r}(\chi_1-\chi_2)\right).
\end{equation}
We now expand the product in (\ref{prodKq2}). The term :
\begin{align*}
\int_{z \in \gamma_\tau^{<,1}}\int_{w \in \gamma_\tau^{>,2}}... \times \int_{z' \in \gamma_\tau^{<,2}}\int_{w' \in \gamma_\tau^{>,1}}...
\end{align*}
is by construction dominated by any polynomial in $r$. The estimates (\ref{ineqres}) and (\ref{ineq2}) imply that the product of the integrals of the residues is smaller than :
\begin{align*}
\frac{Cr^2}{|\chi_1-\chi_2|^2},
\end{align*}
while the combinations of (\ref{ineqres}) and (\ref{ineqint2}), and (\ref{ineq2}) and (\ref{ineqint}) entail that the remaining terms are smaller than :
\begin{align*}
\frac{Cr}{|\chi_1-\chi_2|}.
\end{align*}
Lemma \ref{lemdecK} is proved. $\square$.
\subsection{Proof of Lemma \ref{lemdiag}}
Lemma \ref{lemdiag} is proved using Lemma \ref{lemcv} and proposition \ref{propcont}. Let $K \subset \R^2$ be a compact set, and let $m,m' \subset E$ be finite. Let $(\tau,\chi) \in A\cap K$. By Lemma \ref{lemcv}, we have that :
\begin{multline*}
 \E_r \left[ c_{\frac{1}{r}(\tau,\chi)+m} c_{\frac{1}{r}(\tau,\chi)+m'}\right] - \E_r \left[ c_{\frac{1}{r}(\tau,\chi)+m}\right] \E_r \left[ c_{\frac{1}{r}(\tau,\chi)+m'} \right] \\
 =  \E_{(\tau,\chi)} \left[ c_{m} c_{m'}\right] - \E_{(\tau,\chi)} \left[ c_m\right] \E_{(\tau,\chi)} \left[ c_{m'} \right]+O(r).
\end{multline*} 
Now, by proposition \ref{propcont}, the function :
\begin{align*}
(\tau,\chi) \mapsto  \E_{(\tau,\chi)} \left[ c_{m} c_{m'}\right] - \E_{(\tau,\chi)} \left[ c_m\right] \E_{(\tau,\chi)} \left[ c_{m'} \right]
\end{align*}
is bounded, as long as $(\tau, \chi)$ belong to a compact set. Lemma \ref{lemdiag} is proved, and Theorem \ref{thm1} is completely proved. $\square$.

\end{document}